\documentclass[11pt]{article}

\usepackage{amssymb,amsmath}

\newtheorem{theorem}{Theorem}[section]

\newtheorem{Lemma}[theorem]{Lemma}
\newtheorem{defi}[theorem]{Definition}

\newcommand{\non}{\nonumber}
 \newcommand{\cR}{{\mathbb R}}

\newcommand{\f}{\frac}

\newcommand{\eq}[1]{\mbox{\rm {(\ref{#1})}}}

\numberwithin{equation}{section}

\usepackage{epsfig}

\title{\Large\bf Solvability via viscosity solutions for a model of phase transitions driven by
configurational forces}

\author{\small\sc Peicheng Zhu$^{1,2}$\thanks{ E-mail: zhu@bcamath.org}\ \\
\small $^1$ Basque Center for Applied Mathematics (BCAM)\\
\small Building 500, Bizkaia Technology Park \\
\small E-48160 Derio, \ \
\small Spain\\
\small $^2$ IKERBASQUE, Basque Foundation for Science\\
\small  E-48011 Bilbao,\ \ Spain}

\date{ }

\begin{document}

\maketitle

\begin{abstract}

This article is concerned with an initial boundary value problem for an elliptic-parabolic coupled system
arising in martensitic phase transition theory of elastically deformable solid materials, e.g., steel.
This model was proposed in \cite{Alber04}, and investigated in \cite{Alber06} the existence of
weak solutions which are defined in a standard way, however the key
technique used in \cite{Alber06}  is not applicable to
multi-dimensional problem. One of the motivations of this study is to solve this multi-dimensional problem,
 and another is to investigate the sharp interface limits. Thus we define weak solutions in a way, which
 is different from \cite{Alber06}, by using the notion of viscosity solution. We do prove
successfully the existence of weak solutions in this sense for one dimensional problem,
 yet the multi-dimensional problem is still open.
 %Compared with \cite{Alber06}, the proof is considerably simplified.

\end{abstract}

%%%%%%%%%%%%%%%%%%%%%%%%%%%%%%%%%%%%%%%%%%%%%%%%%%
%%%%%%%%%%%%%%% Section 1 %%%%%%%%%%%%%%%%%%%%%%%%
%
\section{Introduction}
\setcounter{equation}{0}

In this article we shall investigate an initial-boundary
 value problem of a new model which describes martensitic phase transitions in elastically
deformable solid materials, and such phase transitions are driven by configurational forces.
 To formulate this problem, we firstly introduce some notations.
 Let $\Omega$ be an open bounded domain in $\cR^3$ with smooth boundary
$\partial \Omega$. It represents the points of a material body. Define
 $Q_t=(0,t)\times \Omega$. We use unknown functions: $u=u(t,x)$ is the displacement
at time $t$ and position $x$, $T$ is the Cauchy stress tensor, and $S$ is an order parameter which means
that if $S$ takes the values that are approximately equal to $0$ and
$1$, then the material is in two different phases, say $\gamma$ and
$\gamma^\prime$, respectively. Then the system reads
\begin{eqnarray}
 -{\rm div}_x\, T(t,x) &=& b(t,x),
 \label{eq1}  \\
 T(t,x) &=& D\left(\varepsilon(\nabla_x u)-\bar\varepsilon S \right)(t,x),
 \label{eq02} \\
 S_t(t,x) &=& -c\left(-T\cdot\bar\varepsilon + \hat\psi^\prime(S) -\nu
 \Delta_xS\right)|\nabla_xS|(t,x)
 \label{eq2}
\end{eqnarray}
which must be satisfied in $Q_{t}$. We prescribe the following Dirichlet boundary and initial conditions
\begin{eqnarray}
 u|_{[0,t ]\times \partial \Omega}&=& 0,
 \label{eq3a}\\
 S|_{[0,t ]\times \partial \Omega}&=& 0,\\
 S|_{\{t=0\}\times {\bar \Omega} } &=& S_0.
 \label{eq3}
\end{eqnarray}

In this model, $c,\ \nu$ are positive constants, $D$ is  the linear, positive definite symmetric elasticity tensor.
 We have chosen the free energy $\psi = \psi(\varepsilon,S,\nabla_xS)$   given by
\begin{eqnarray}
 \psi(\varepsilon,S,\nabla_xS)=\f12 D(\varepsilon-\bar\varepsilon S)
 \cdot (\varepsilon-\bar\varepsilon S ) +\hat\psi(S)+\frac\nu2|\nabla_xS|^2,
 \label{freeenergy}
\end{eqnarray}
 and $\psi_S$ is the derivative, with respect to $S$, of $\psi$. The scalar product of two matrices
 $\sigma,\ \tau$ is denoted by $\sigma\cdot \tau = \sum_{i,j=1}^3\sigma_{ij} \tau_{ij}$. There holds
  the relation $\psi_S(\varepsilon, S)=-T\cdot \bar\varepsilon + \hat\psi^\prime(S)$.
  $\varepsilon$ is the strain tensor defined by
$\varepsilon=\varepsilon(\nabla_x u)=\frac12\left(\nabla_x u+\, ^t(\nabla_x u) \right)$,
 and the upper-script~$^t$ denotes the transpose of a matrix. $ \bar\varepsilon  $  is called
the misfit strain. The function $\hat\psi(S)$ is chosen as a double-well
potential for which we assume that
\begin{eqnarray}
 & & \hat\psi(S)\  {\rm is\ smooth\ and\ has\ two\ minima\ at\ }  S=0
 {\rm\ and\ }  S=1, \  {\rm and}\non\\
 & & {\rm one
 \ maximum\ at }\ \hat S  \ {\rm between}\ 0\ {\rm and} \ 1,
 \label{assumppsi}\\
 & &\hat\psi'(S)>0,\ {\rm if}\ S\in (0,\hat S )\cup (1,\infty);
 \ \hat\psi'(S)<0,\ {\rm if}\ S\in (\hat S ,1)\cup (-\infty,0).
 \non
\end{eqnarray}
Finally, $b = b(t,x) $ is a given volume force.

\medskip
This model was formulated in \cite{Alber04} by employing the second law of thermodynamics
and a formula (see e.g. \cite{AK90,Gurtin,Maugin}) of configurational forces. Our model differs from the celebrated Allen-Cahn model
 (which is also called Ginzburg-Landau) by the gradient term $|\nabla_xS|$. The reason is that in the Allen-Cahn model, the driving force for
the motion of interfaces is the mean curvature, while the motion of interfaces considered in this
paper is driven by configurational forces.  We mention the key ideas of the derivation. There are two main
types of phase transition models: sharp interface model and phase field model. Our model is derived from a
 sharp interface model: Assuming that the jump of $S$, across the interface of two phases, becomes smaller and smaller,
 we see that the equation governing the interface approaches to a Hamilton-Jacobi equation
 $S_t= -c\, \psi_S |\nabla_x S|$ which is
 a fully nonlinear equation, thus is difficult to deal with and its solution may develop singularities.
 A usual way for regularizing it is to add an artificial term (for instance, $\nu\Delta_x S$)
  as in the theory of conservation laws, but this
technique does not work in our case. We then think of another type of models, i.e. phase field model, to regularize such an
equation. To formulate a phase field model, we choose the free energy \eq{freeenergy}, and also need a suitable flux which
 can be chosen  in the form
\begin{eqnarray}
 q = q(u_t,T,\nabla_xS,S_t)=T\cdot u_t +\nu S_t\nabla_x S.
 \label{flux}
\end{eqnarray}
Then by straightforward computations, we see that if the equations (\ref{eq1}) -- (\ref{eq2}) are satisfied, then the following
Clausius-Duhem inequality is satisfied
\begin{eqnarray}
 \frac{d}{dt}\, \psi(\varepsilon,S,\nabla_xS)-{\rm div}_x\, q - b\cdot u_t\le 0.
 \label{2ndlaw}
\end{eqnarray}
Hence, we assert the validity of the second law of thermodynamics.
For the details of the formulation of this model, we refer to the appendix of this paper, or the articles
\cite{Alber00,Alber06,Alber04}.

\medskip
The aim of this article is to propose a suitable concept of weak solutions that
works for multi-dimensional problem and that makes the
investigation of sharp interface limit (as $\nu$ goes to $0$) easier,
% to deal with,
 then to prove the existence of such defined weak solutions for problem \eq{eq1} -- \eq{eq3}.
 There are two most well-known concepts of weak solutions to partial
differential equations: The first one is the notion of usual weak solutions that are
defined by employing test functions and the technique of integration
by parts, and the second one is the conception of viscosity solutions
developed by Crandall and Lions in 1983, see \cite{CL83}, etc. In
this article, we define a weak solution by combining these two notions of weak
solutions. To understand why we need two concepts of weak solutions, we first
 investigate the features of this model. Our model consists of a subsystem of linear elasticity
 and a partial differential equation that is degenerate and has  strong nonlinearity and non-smooth coefficients. The
 one space dimensional initial-boundary value problem for this model has been studied
 in \cite{Alber06}, in which we define a weak solution
in a usual way by using a simple technique that makes us possible
to rewrite the principle part of the equation of the order parameter in
a divergence form, i.e. $\nu S_{xx}|S_{x}| = \frac\nu2(S_{x}|S_{x}|)_x$.

However such a technique fails for the corresponding
 multi-dimensional problem of this model, namely $\nu \Delta_x S |\nabla_{x}S|$ can't be rewritten in
 a divergence form. Thus the notion of usual weak solutions is not suitable for
 this problem because we can not reduce the order of weak derivatives of
 the unknown by integration by parts. This is one of the difficulties in solving our
 model. Another one is that the maximum principle, which plays a
 crucial role in the theory of viscosity solutions, is not valid for the whole system of
 equations considered here. So it is not suitable to define weak solutions by using the notion of
 viscosity solutions only. Therefore one of two purposes of this article is to propose a suitable notion of
 weak solutions to this multi-dimensional problem.
 %, though we can only solve one dimensional  problem up to now.
 The second purpose is that we shall use our new notion of weak solutions to study,
 in the future,  a very interesting problem, i.e. the sharp interface limit of our model. Such a problem however
  may be difficult   under the framework of the standard weak solution, since the sharp interface problem
  has a fully nonlinear equation of the order parameter.

The above consideration leads us to propose a suitable notion of
generalized solutions to our system by using both notions of weak solutions:
 we define weak solutions in the usual sense for
the subsystem of elasticity, and use viscosity solutions to define weak solutions to the order
parameter equation.  Then we construct a sequence of  solutions to an  approximate initial boundary value problem
of the system. Applying some compactness lemma we can show that the limit of the approximate solutions
is just weak solutions in our sense. Though only the one space dimensional problem is
solved up to now, we believe this technique works for the multi-dimensional case too.
The other interesting open problems in this field include: The sharp interface limit of our model, and the
relationship between weak solutions defined in
this article and the ones  in \cite{Alber06}, respectively.

\medskip
We are now going to study the definition and existence of weak solutions in a suitable sense
 to problem (\ref{eq1}) -- (\ref{eq3}) in one space dimension, though the definition and
some {\it a-priori} estimates are still valid for multi-dimensional problem. We shall see
 later on that the proof of the existence of weak solutions in this article is significantly
 simpler than that in \cite{Alber06}.

\noindent {\bf Statement of the main result.} From now on we
 assume that all functions only depend on the variables $x_1$ and $t$, and,
to simplify the notation, denote $x_1$ by $x$. The set $\Omega=(a,d)$ is a
bounded open interval with constants $a<d$. We write $Q_{t_e}:=(0,t_e) \times
\Omega$, where $t_e$ is a positive constant, and define
$$
 (v,\varphi)_Z = \int_Z v(y)\varphi(y)\, dy\,,
$$
for $Z = \Omega$ or $Z = Q_{t_e}$. If $v$ is a function defined on
$Q_{t_e}$ we denote the mapping $x \mapsto v(t,x)$ by $v(t)$. If
no confusion is possible we sometimes drop the argument $t$ and
write $v = v(t)$. We still allow that the material points can be
displaced in three directions, hence $u(t,x)\in\cR^3$, $T(t,x)\in
{\cal S}^3$ and $S(t,x)\in\cR$, where ${\cal S}^3$ is the set of $3\times 3$ symmetric
matrices. If we denote the first column of the matrix
$T(t,x)$ by $T^1(t,x)$ and set
$$
 \varepsilon(u_x)=\frac12\left((u_x,0,0)+\,^t(u_x,0,0)\right)\in {\cal S}^3,
$$
then with these definitions the equations (\ref{eq1}) -- (\ref{eq2}) in
the case of one space dimension can be written in the form
\begin{eqnarray}
 -T^1_{x} &=& b,
 \label{eq1a}\\
 T &=& D(\varepsilon(u_x)-\bar\varepsilon S),
 \label{eq1aa} \\
 S_t &=& c\left(T\cdot\bar\varepsilon-\hat\psi^\prime(S) + \nu S_{xx}\right)|S_x|,
 \label{eq2a}
\end{eqnarray}
which must be satisfied in $ Q_{t_e}$.
%Since the equations (\ref{eq1a}), (\ref{eq1aa}) are linear, the
%inhomogeneous Dirichlet boundary condition for $u$ can be reduced in the
%standard way to the homogeneous condition. For simplicity we thus assume that
%$\gamma=0$.
 The boundary and initial conditions therefore are
\begin{eqnarray}
 u(t,x)=0, && (t,x)\in [0,t_e]\times\partial\Omega,
 \label{eq3aa}\\
 S(t,x)=0, && (t,x)\in [0,t_e]\times\partial\Omega, \\
 S(0,x)=S_0(x), && x\in  \Omega.
 \label{eq4a}
\end{eqnarray}

\medskip
To define weak solutions to problem (\ref{eq1a}) -- (\ref{eq4a}),
we first introduce some definitions on semi-continuous functions.
Let $f=f(x)$ be a real function defined in $U\subset \cR^N$ with $N\in\{1,2,5\}$.
We denote the
so-called upper semi-continuous envelope of $f$ by
\begin{equation}
  f^*(x): U\to \cR\cup \{-\infty,+\infty\}
\end{equation}
which is defined by
\begin{equation}
  f^*(x):=\lim_{r\downarrow 0}{{\rm ess}\sup}_y\{f(y)\mid y\in U,
 \ \ |x-y|\le r\}.
\end{equation}
Obviously, $f^*(x)$ is upper semi-continuous. And $f_*(x):=-(-f)^*(x)$
is called lower semi-continuous envelope of $f$.

We define the Hamiltonian $H_T$ which depends on the unknown $T$ by
\begin{equation}
 H_T(t,x,p,q,r)=c\left( T(t,x)\cdot \bar \varepsilon-\hat\psi^\prime(p)+\nu\, r\right)|q|,
 \label{H1}
\end{equation}
where, $(t,x)\in Q_{t_e},\, p,\, q, \, r\in \cR$, so $(t,x,p,q,r)\in\cR^5$.
It is easy to show that if $T$ is a continuous function in $(t,x)$
 and $\hat\psi^\prime$ is continuous in $S$,
 then we have that $H_T$ is continuous in $(t,x,p,q,r)$, thus
\begin{equation}
 (H_T)^*(t,x,p,q,r)=(H_T)_*(t,x,p,q,r)=H_T(t,x,p,q,r).
 \label{H}
\end{equation}

\medskip
We now can introduce the notion of weak solutions for our problem. In what follows we
shall assume that  $p$ is a real number such that
\begin{equation}
 p>1.
 \label{P}
\end{equation}
%\vskip0.2cm
\begin{defi}\label{D1.1}
A function  $(u,T,S)$ which satisfies that
\begin{equation}
 (u,T,S)\in L^\infty(0,t_e;H^{1 }_0(\Omega))
 \times L^\infty(0,t_e;L^2(\Omega))\times L^\infty(\bar Q_{t_e}),
 \label{Properties}
\end{equation}
is called a weak solution to system (\ref{eq1a}) -- (\ref{eq4a}) if
\begin{description}
\item \quad I)  for almost every  $t\in [0,t_e]$, equations \eq{eq1a}, \eq{eq1aa}
and the boundary condition \eq{eq3aa} are  satisfied  weakly.
%\begin{equation}
% \int_\Omega T^1(t,x) \cdot \varepsilon (v_x)(x)dx=\int_\Omega b(t,x)\cdot v(x) dx,
%  \quad \int_\Omega (T(t,x) -  D\left(\varepsilon (u_x )-\bar\varepsilon S \right))(t,x)
%  w (x)dx =0
%\end{equation}
%for any function $v  \in H_0^1(\Omega,\cR^3),\ w \in L^2(\Omega)$.
%, and relation (\ref{eq1aa}) is satisfied in the sense of distribution.
 %
\item \quad II) $S$ is a viscosity solution to equation (\ref{eq2a}), if  $S$
 satisfies both i) and ii) below:

\indent \quad\quad  i)  $S$ is a sub-viscosity solution to equation (\ref{eq2a}),
 i.e.  for any  function $\phi(t,x)$ in $C^{1,2}(\bar Q_{t_e})$,
if $S^*-\phi$ attains its local maximum at $(\tau, y)$, then
\begin{equation}
 \phi_t(\tau, y)\le (H_T)_*(\tau, y,S^*(\tau, y),\phi_x(\tau, y),\phi_{xx}(\tau, y)),
\end{equation}
and $S^*(0,x)\le S_0(x)$;

\indent \quad\quad  ii)  $S$ is a super-viscosity solution to
Eq. (\ref{eq2}), i.e. for any
 function $\phi(t,x)$ in
$C^{1,2}(\bar Q_{t_e})$, if $S_*-\phi$ attains its local minimum at
$(\tau, y)$, then
\begin{equation}
 \phi_t(\tau, y)\ge (H_T)^*(\tau, y,S_*(\tau, y),\phi_x(\tau, y),\phi_{xx}(\tau, y) ),
\end{equation}
and $S_*(0,x)\ge S_0(x)$.
\end{description}

\end{defi}

Now we are able to state our main result  as follows.
%%%%%%%%%%%%%%%%%   Main Theorem  %%%%%%%%%%%%%%%%%%%%
\begin{theorem}  \label{T1.2}
Suppose that $b,\ b_t\in C([0,t_e];L^2(\Omega))$ for any given
positive constant $t_e$, and that $S_0\in  H^1_0(\Omega)$.
% satisfies
%$0\le S_0(x)\le 1$ for all $x\in \bar\Omega$, and
%the compatibility condition  $S_0|_{\partial \Omega}=0$.
Furthermore, we assume that the function $\hat\psi$ satisfies the
assumption (\ref{assumppsi}).

Then there exists a weak solution $(u,T,S)$ to problem (\ref{eq1a})
-- (\ref{eq4a}) in the sense of Definition~\ref{D1.1}, and in
addition to (\ref{Properties}), we have that the solution satisfies
$$
 S\in C(\bar Q_{t_e}).
$$
%and $0\le S(t,x)\le 1$ for all $(t,x)\in \bar Q_{t_e}$.

\end{theorem}

\bigskip
Our notion of generalized solutions is a combination of the concept of usual weak
solutions and the notion of viscosity solutions.
%(applied to the linear second order elliptic subsystem of
%the unknown $u$ in our system)
% (applied to the nonlinear scalar part, i.e.
%for the equation of order parameter $S$).
This idea comes partly from some discussions with Prof. Alber and partly from the paper
by Giga, Goto and Ishii \cite{GGI92} which is concerned with the global
existence of weak solutions, however without uniqueness, to the system consisting of
a semi-linear diffusion equation in two disjoint open sub-domains denoted by
$\Omega_\pm(t)$ of one simply connected domain $\Omega$ (The complement of
 union of these two parts is so-called the interface $\Gamma(t)$), and a nonlinear
 interface equation. The system is composed of the interface equation
\begin{equation}
 V=W(v)- c\, \kappa, \ \ {\rm on\ \ } \Gamma(t)
 \label{Vformula}
\end{equation}
and the diffusion equations
$$
v_t=\nu \Delta v+g_\pm(v), \ \ {\rm for }\ \ x\in \Omega_\pm(t), \ \ t>0.
$$
Here, $V=V(t,x)$ is the speed of $\Gamma(t)$ at $x\in \Gamma(t)$ in the normal
 direction of $n$ from $\Omega_+(t) $ to $\Omega_-(t) $. $\kappa$ is the mean
 curvature of $\Gamma(t)$ at $x\in \Gamma(t)$, $v$ is the density.
 And $c, \nu $ are positive
constants, $ W, g_\pm $ are given bounded continuous functions over
$\cR$. Note that in the work \cite{GGI92}, the driving force for
the motion of an interface is due to the mean curvature (see formula \eq{Vformula}),
 while the motion of an interface considered in this article
 is driven by configurational forces and the motion is
 governed by $ V[S] = cn\cdot [E]n$ (the sharp interface case),
 where $E$ is the Eshelby tensor, an energy-momentum tensor,
 see  \cite[pp. 753-767]{Eshelby}.

\bigskip
We recall the literature related to our results.
There have been many papers on the theory of viscosity solutions
 since the notion of viscosity solution was proposed in 1983 by Crandall
 and Lions~\cite{CL83}. This notion is applicable to the scalar equations
 or the weakly coupled systems, for which the maximum principle  holds. Hence, the
comparison theorem is valid, this plays an
 important role in the proof of uniqueness of viscosity solution.
 For an overview of the theory, we refer for instance to Capuzzo Dolcetta and Lions \cite{CDL97},
Ishii and Lions \cite{IL90}, Crandall, Ishii and Lions \cite{CIL92},
Jensen\cite{Jensen88}, Crandall and Lions~\cite{CL86},
Ishii \cite{Ishii87}, Souganidis \cite{Souganidis85} for the scalar
equation case,  and to Engler and Lenhart \cite{EL91}, Ishii and
Koike \cite{IK91}, etc. for the system case,   and the references are
cited therein. For the background of our model and mathematical results related this article,
we refer the reader to work by Alber and/or Zhu \cite{Alber00,Alber06,Alber04,Alber07,Alber09,Alber10,Zhu10},
Kawashima and Zhu \cite{kawa}, Ou and Zhu\cite{Ou10}.

The main difficulties and our strategies in the proof of
Theorem~\ref{T1.2} are as follows: {Firstly}, the definition of weak
solutions is a new problem since our system comprises of a linear
elliptic system of $u$ and a nonlinear equation of $S$ which can not
be rewritten in the divergence form.
{Secondly}, the equation for the order parameter is degenerate and its
coefficients is not smooth.
 To overcome these difficulties, we make a suitable smooth approximation of the non-smooth term
which leads the equation of the  order parameter to a uniformly parabolic equation with smooth
 coefficients. We employ the energy estimates to discuss the limits of approximate solutions.

\vskip0.2cm
The remaining of this article is organized as follows. In Section~2 we state an approximate
 initial boundary value problem, and apply the existence theorem in the book by Ladyzenskaya et al.
\cite{Ladyzenskaya} to prove existence of classical solution to this approximate problem.
Then we derive in Section~3 the uniform {\it a priori} estimates which are independent of a small parameter
$\kappa$ for the approximate solutions. Then we apply the {\it a priori} estimates, a lemma of the Aubin-Lions type  and a
theorem on the stability of viscosity solutions to discuss the limits and prove the
existence of weak solutions in the sense of Definition~1.1. Finally Section~4 we present briefly in the appendix
 the derivation of our model.

%%%%%%%%%%%%%%%%%%%%%%%%%%%%%%%%%%%%%%%%%%%%%%%%%%
%%%%%%%%%%%%%%% Section 2 %%%%%%%%%%%%%%%%%%%%%%%%
%
\section{Existence of solutions to the modified problem}
\setcounter{equation}{0}

In this section, we are going to study an approximate initial-boundary value problem and
show that it has a classical solution for any fixed positive constant $\kappa$. Since we shall
let $\kappa$ go to zero, we may assume, without loss of generality, that
$$
 0<\kappa<1.
$$
Let
$
\chi\in C_0^\infty(\cR^2,[0,\infty))
$
be a function satisfying
$
\int_{-\infty}^{\infty} \chi(t,x)dtdx=1.
$
We set
$$
\chi_\kappa(t,x):=\frac1{\kappa^2} \chi\left(\frac{t}{\kappa},\frac{x}{\kappa} \right),
$$
and for $b\in L^\infty(Q_{t_e},\cR)$ we define
\begin{eqnarray}
 (\chi_\kappa*b)(t,x)=\int_{0}^{t_e} \chi_\kappa(t-s,x-y)b(s,y)ds dy.
 \label{convolution}
\end{eqnarray}
We smooth the term $|S_x|$  as follows
\begin{eqnarray}
 | S_x|_\kappa=\sqrt{ |S_x|^2+\kappa^2 },
 \label{smoothed}
\end{eqnarray}
and choose a sequence $ S_0^\kappa$ such that
\begin{eqnarray}
 S_0^\kappa\in C_0^\infty(\Omega),\quad \|S_0^\kappa - S_0\|_{H^1(\Omega)}\to 0
\label{initial}
\end{eqnarray}
as $\kappa\to 0$ since $ C_0^\infty(\Omega)$ is dense in $H^1_0(\Omega)$.

Then the smoothed initial boundary value problem of
(\ref{eq1a}) -- (\ref{eq4a}) turns out to be
\begin{eqnarray}
 -T^1_{x} &=& \chi_\kappa*b,
 \label{eq1appro1}  \\
 T&=&D\left(\varepsilon( u_x)-\bar\varepsilon S \right),
  \\
 S_t &=& c\, \nu |S_x|_\kappa  S_{xx}
 +c\left(T\cdot \bar\varepsilon-\hat \psi^\prime(S)\right)(|S_x|_\kappa - \kappa).
 \label{eq2appro1}
\end{eqnarray}
and the boundary and initial conditions become
\begin{eqnarray}
 u|_{[0,t_e]\times\partial\Omega}&=&0,\\
 S|_{[0,t_e]\times\partial\Omega}&=&0,
 \label{BDappro1}\\
 S|_{\{0\}\times\bar\Omega}&=&S_0^\kappa.
 \label{IDappro1}
\end{eqnarray}
By the choice of $S_0^\kappa$, we see that the compatibility
condition $S_0^\kappa|_{\partial\Omega}=0$ is met.

\medskip
\noindent{\bf Remark 2.1.} {\it There are some other ways, which are different from \eq{smoothed}, to smooth
the function $|p|$. We need only to require that the smoothed  equation \eq{eq2appro1} for the order parameter
meets the assumptions of the maximum principle.
}

\medskip
To prove the existence of classical solution to the
approximate problem (\ref{eq1appro1}) -- (\ref{IDappro1}), we
employ the Leray-Schauder  fixed-point theorem (see, e.g. \cite{Ladyzenskaya}) to the following problem
\begin{eqnarray}
 -T^1_{x} &=& \lambda\,\chi_\kappa*b,
 \label{eq1appro}  \\
 T &=& D\left(\varepsilon( u_x)-\lambda\, \bar\varepsilon S \right),
 \label{eq1appro11}\\
 S_t &=& c\, \nu |S_x|_\kappa  S_{xx}
 +c\left(T\cdot \bar\varepsilon-\hat \psi^\prime(S)\right)(|S_x|_\kappa - \kappa).
 \label{eq2appro}
\end{eqnarray}
Here $\lambda \in [0,1]$. The boundary and initial conditions are
\begin{eqnarray}
 u|_{[0,t_e]\times\partial\Omega}&=&0,\\
 S|_{[0,t_e]\times\partial\Omega}&=&0,
 \label{BDappro}\\
 S|_{\{0\}\times\bar\Omega}&=&\lambda\, S_0^\kappa.
 \label{IDappro}
\end{eqnarray}

Define for any $\hat S\in {\cal B}:= C^{1+\frac{\alpha}2, 2+\alpha}(\bar Q_{t_e})$
(here $0<\alpha<1$) a mapping $P_\lambda : [0,1]\times {\cal B}\to {\cal B};\quad \hat S\mapsto S$
where $S$ is a solution obtained by the following procedure:

i) For any fixed $\hat S$, it is easy to find a unique
 solution $(u,T)$ which depends on $\hat S$, to the following boundary value problem for almost
 every given $t$
\begin{eqnarray}
 - T^1_{x} &=& \lambda\,\chi_\kappa*b,
 \label{eq1'}  \\
 T&=&D\left(\varepsilon( u_x)-\lambda\, \bar\varepsilon  \hat S\right),\\
 u|_{\partial\Omega}&=&0.
 \label{eq1a'}
\end{eqnarray}

ii) Then inserting this $T$ into equation (\ref{eq2appro}) we can
 obtain a unique classical solution $S$ to problem (\ref{eq2appro}),
(\ref{BDappro}) and (\ref{IDappro}).

\medskip
With the help of some a priori estimates,  we then obtain
\begin{theorem} \label{T2.1}
Suppose that all the assumptions in Theorem~\ref{T1.2} are met, and
the compatibility conditions $S_{0}=S_{0,x}=S_{0,xx}=0$ at $x=a,d$
are satisfied.

Then for any fixed $\kappa>0$, there exists a unique classical solution $(u,T,S)\in C^{2,1}(\bar{Q}_{t_e})
\times C^{1,1}(\bar{Q}_{t_e}) \times C^{2+\alpha,1+\alpha/2}(\bar{Q}_{t_e})$
to  problem (\ref{eq1appro}) -- (\ref{IDappro}) which satisfies
\begin{eqnarray}
 S_{tx}\in L^2(Q_{t_e}).
 \label{Sxt}
\end{eqnarray}

\end{theorem}

\medskip
\noindent{\bf Remark 2.2.} {\it The compatibility conditions in
Theorem~2.1 are  different from usual ones and they are
derived as follows: From the system and initial data, there must
hold
\begin{eqnarray}
 T (0,x)|_{x=a,d} - \left. D \varepsilon( u_x(0,x))
 \right|_{x=a,d} & = & 0,\non\\
 \ \label{compty}\\
 \left. \nu |S_{0,x}|_\kappa  S_{0,xx} +
 T(0,x)\cdot \bar\varepsilon (|S_{0,x}|_\kappa - \kappa)
 \right|_{x=a,d} &=& 0.\non
\end{eqnarray}
Note that the values of $u_x(0,x)$ at boundary can be arbitrary, so
is  $T (0,x)$. Thus from  the definition of the function  $|\cdot|_\kappa$ we see
 that the second term of (\ref{compty}) is satisfied provided that $S_{0,x}=S_{0,xx}=0$ at $x=a,d$.

}

\bigskip
The proof of Theorem~\ref{T2.1} is equivalent to  that $P_\lambda$ has a fixed point for $\lambda=1$. 
It is clear that if $\lambda=0$, $P_\lambda$ maps any $\hat S\in {\cal B}$ to $S\equiv 0$, i.e. 
$P_0 \hat S =0$. Thus by the Leray-Schauder fixed-point theorem, we see it remains to derive the  
estimates that are stated in Lemma~2.2 -- Lemma~2.5 and Lemma~2.7, from which we conclude the compactness
of $P_\lambda$. To obtain those a priori estimates, we assume that there exists a  classical solution 
$(u,T,S)\in C^{2,1}(\bar{Q}_{t_e}) \times C^{1,1}(\bar{Q}_{t_e}) \times C^{2+\alpha,1+\alpha/2}(\bar{Q}_{t_e})$
to problem (\ref{eq1appro}) -- (\ref{IDappro}) such that $S_{tx}\in L^2(Q_{t_e})$.

\medskip
Firstly, writing $T$ in terms of $S$ we reduce system (\ref{eq1appro}) -- (\ref{eq2appro}) 
into a scalar equation with a non-local term, and  apply the maximum
principle this equation to obtain
\begin{Lemma} \label{L2.2}
There holds for $t_e>0$
\begin{eqnarray}
 \|S \|_{L^\infty(Q_{t_e})}&\le& \bar C.
 \label{est2}
\end{eqnarray}

\end{Lemma}

In this lemma and in the follows context, we denote by $\bar C$ a constant which is
 independent of $\kappa$,  but may depend on $\nu$, while a constant $C$   
 may depend on both $\kappa$ and $\nu$.

\medskip
\noindent{\it Proof.} To make use of the maximum principle, we solve $(u,T)$ in terms of
$S$ from the first two equations \eq{eq1appro} -- \eq{eq1appro11}, provided that $S$ is given. 
Then the whole system can be reduced into a single equation,
but with a nonlocal term. We need some notations as
used in \cite{Alber04}. Let $\hat{\cal S}^3$ be the subspace of all matrices $A\in {\cal S}^3$
with $A_{ij}=0$ for $i,j=2,3$. The orthogonal space to $\hat{\cal S}^3$ is denoted by
$\tilde{\cal S}^3$. It consists of $A\in {\cal S}^3$ satisfying $A_{i1}=A_{1i}=0$ for all
$i=1,2,3$. Note that $\varepsilon(u_x)\in \hat{\cal S}^3$. Let $\hat P$ be
 the canonical projection of $ {\cal S}^3$
 into $\hat{\cal S}^3$. Since $D : {\cal S}^3\to {\cal S}^3$ is a positive definite linear mapping,
 $\langle \sigma,\tau\rangle = D\sigma\cdot \tau$ defines a scalar product on $ {\cal S}^3$. The
 projection of $ {\cal S}^3$ onto $ \hat{\cal S}^3$, which is orthogonal with respect to this
 scalar product is denoted by $\hat Q$. These definitions imply that
$$
 {\rm ker }\, \hat Q =D^{-1}\tilde{\cal S}^3 = D^{-1} {\rm ker }\, \hat P.
$$
Define further that
$$
 \varepsilon^* = \hat Q\bar \varepsilon,\ u^* = (\varepsilon^*_{11}, 2\varepsilon^*_{21}, 2\varepsilon^*_{31}),
$$
we then obtain
\begin{eqnarray}
 u(t,x) &=&  \lambda\,u^*\left(\int_a^x S(t,y)dy - \frac{x-a}{d-a}\int_a^d S(t,y)dy\right) + w(t,x),\label{uformula}\\
 T(t,x) &=& D(\varepsilon^* - \bar\varepsilon) \lambda\,S(t,x) -  \frac{D \varepsilon^*}{d-a}\int_a^d \lambda\,S(t,y)dy  + \sigma(t,x),
 \label{Tformula}
\end{eqnarray}
where  the function $(w(t,\cdot),\sigma(t,\cdot))$ (here $t$ is temporarily regarded as a parameter) is the
 unique solution of the following boundary value problem
\begin{eqnarray}
 - \sigma_{1x}(x) &=& \hat b(x)\ \ {\rm in}\ \ \Omega, \\
 \sigma (x) &=& D\varepsilon(w_x(x))\ \ {\rm in}\ \ \Omega, \\
 w(a) &=& \hat f(a), \ w(d)  =  \hat f(d)
\end{eqnarray}
and $\hat b =  \lambda\, b (t)$, $\hat f  \equiv 0$. Note that $u^*\in \cR^3$, $\varepsilon^*\in {\cal S}^3$
depend only on the misfit strain $\bar\varepsilon$. Inserting the formula  of $u,\ T$ into equation \eq{eq2appro} yields that system
 \eq{eq1appro} --  \eq{eq2appro} is reduced into a single equation for $S$ with a nonlocal term. Invoking the
 definition of $|p|_\kappa$, we see that the assumptions required by the maximum principle are satisfied.
Thus we can apply the maximum principle to this single equation and the proof of this lemma is complete.

\medskip
Next we can derive the following estimates for the derivatives of $S$.
\begin{Lemma}
\label{L2.6}
There holds  for any $t\in [0,t_e]$ that
\begin{eqnarray}
 \|S_x(t)\|^2 + \int_0^t\int_\Omega | S_x|_\kappa  |S_{xx}|^2dx d \tau &\le& \bar C,
 \label{es2a}\\
 \int_0^t\int_\Omega\left(\left(| S_x |_\kappa |S_{xx} |\right)^\frac43 + |S_t|^\frac43\right) dx d \tau &\le& \bar C.
 \label{es2b}
\end{eqnarray}
Here and hereafter, we denote the $L^2$-norm over $\Omega$ by $\|\cdot\|$.

\end{Lemma}

\noindent{\it Proof.} By definition we have the property $|p|_\kappa
\ge \kappa$, from which we obtain
$$
 0\le |p|_\kappa- \kappa \le |p|_\kappa + \kappa\le 2|p|_\kappa.
$$
Using estimate \eq{est2} and formula \eq{Tformula}, recalling the assumptions on $b$,  one concludes that
\begin{eqnarray}
 \|T \|_{L^\infty(Q_{t_e})}\le \bar C.
 \label{Tbound}
\end{eqnarray}
Note that $S_{tx}\in L^2(Q_{t_e})$, for any fixed $\kappa$,  implies that
$$
 \frac12\frac{d}{dt}\| S_x \|^2 = ( S_x,  S_{xt}).
$$
 Multiplying (\ref{eq2appro}) by $-S_{xx} $ and integrating
 the resulting equation with respect to $x$, using integration
by parts, and invoking  the estimates (\ref{est2}) and  (\ref{Tbound})  we get
\begin{eqnarray}
 \frac12\frac{d}{dt}\| S_x \|^2+
 c\, \nu\left( |S_x |_\kappa  S_{xx} ,S _{xx}\right)
 &=& c\left(\left(T \cdot \bar\varepsilon
 -\hat \psi^\prime(S )\right)(|S_x |_\kappa - \kappa),
 -S _{xx}\right) \non\\
 &\le& \bar C\left(  |S_x |^\frac12_\kappa ,\, |S_x |^\frac12_\kappa
 |S _{xx}| \right),
 \label{tempo}
\end{eqnarray}
where we used the notation $(f,g)=\int_\Omega f(x)g(x)dx$.
Applying the Cauchy-Schwarz inequality,
% $|\int_\Omega fgh\, dx|\le  \|f\|_{L^4(\Omega)}  \|g\|_{L^4(\Omega) } \| h\| $,
 we infer from \eq{tempo} that
\begin{eqnarray}
 \frac12\frac{d}{dt}\| S_x \|^2+
 c\, \nu\left( |S_x |_\kappa S_{xx} ,S _{xx}\right)
 &\le & \bar C \|\, |S_x |^\frac12_\kappa\|_{L^1(\Omega)}
 \|\, |S_x |_\kappa^\frac12  S _{xx}\| \non\\
 &\le & \bar C(\|S_x \|^\frac12  + 1)
   \|\,|S_x |_\kappa^\frac12 S _{xx}\|
 \label{Remark31}
\end{eqnarray}
By the Young inequality and the property that $|p|_\kappa\le |p| + \kappa$ for any $\kappa\ge 0$, we derive from \eq{Remark31} that
\begin{eqnarray}
\frac12\frac{d}{dt}\| S_x \|^2+
 c\, \nu\left( |S_x |_\kappa S_{xx} ,S _{xx}\right)
 &\le & \frac{c\, \nu}{2}\|\, |S_x |_\kappa^\frac12S _{xx}\|^2
 +\bar C_\nu(\|S_x \| + 1) \non\\
 &\le & \frac{c\, \nu}{2}\int_\Omega |S _x|_\kappa |S _{xx} |^2dx
 + \bar C_\nu\|S_x \|^2+\bar C.\quad
 \label{Remark31z}
\end{eqnarray}
 Thus we arrive at
\begin{eqnarray}
 \frac{d}{dt}\|S_x \|^2+ c\, \nu\int_\Omega |S_x |_\kappa | S _{xx}|^2dx
 \le \bar C_\nu\|S_x \|^2 + \bar C.
 \label{est2c1}
\end{eqnarray}
Using the Gronwall inequality to (\ref{est2c1}) one can easily obtain (\ref{es2a}).

By the interpolation technique and (\ref{es2a}), we have
that for some $2 > p \ge 1, q=\frac2p$ and $\frac1q +
\frac{1}{q^\prime}=1 $ that
\begin{eqnarray}
& & \int_0^t\int_\Omega\left(|S_{x} |_\kappa|S_{xx} |
  \right)^pdx d\tau \non\\
&= & \int_0^t\int_\Omega \left(| S_{x} |_\kappa\right)^\frac{p}{2}
  \left( \left(|S_{x} |_\kappa\right)^\frac{p}{2}
  |S_{xx} |^p\right) dx d\tau \non \\
&\le&  \left(\int_0^t\int_\Omega
  \left(| S_{x} |_\kappa \right)^\frac{pq^\prime}{2}  dxd\tau
  \right)^\frac{1}{q^\prime}\left(\int_0^t\int_\Omega
  \left(|S_{x} |_\kappa \right)^\frac{pq}{2}
  |S_{xx} |^{pq}  dxd\tau \right)^\frac{1}{q} \non\\
&\le &  \left(\int_0^t \int_\Omega
  \left(|S_{x} |_\kappa\right)^\frac{p}{2-p}
  dx d\tau\right)^\frac{2-p}{2}
  \left(\int_0^t \int_\Omega |S_{x} |_\kappa|S_{xx} |^2
  dx d\tau \right)^\frac{p}{2}.
\label{6}
\end{eqnarray}

Invoking the property that $|p|_\kappa\le |p| + \kappa$ and
inequality (\ref{es2a}) yield that for $\frac{p}{2-p}\le 2$, i.e. $p\le
 \frac43$, the right hand side of (\ref{6}) is bounded.

 Making use of (\ref{6}) (with $p =
 \frac43$) and equation (\ref{eq2appro})  we have
 for  any test function $\varphi\in L^4({Q_{t_e}})$
\begin{eqnarray}
 && \left|(S_t ,\varphi)_{Q_{t_e}} \right|
 = c \left|\left(\nu |S_x |_\kappa
  S_{xx} +  (T \cdot \bar\varepsilon-\hat \psi^\prime(S ) )(|S_x |_\kappa - \kappa),
 \varphi\right)_{Q_{t_e}}\right|\non\\
 & \le & \bar C \left\| |S_x |_\kappa
 S_{xx} \right\|_{L^\frac43(Q_{t_e})} \|\varphi\|_{L^4(Q_{t_e})}\non\\
 && +\bar C \left\| (T \cdot \bar\varepsilon-\hat \psi^\prime(S )  )
 \right\|_{L^4(Q_{t_e})} (\|S_x \|+1)  \|\varphi\|_{L^4(Q_{t_e})}\non\\
 & \le & \bar C\left(\left\| |S_x |_\kappa  S_{xx}
 \right\|_{L^\frac43(Q_{t_e})}+
 \|S_x \| + 1 \right)\|\varphi\|_{L^4(Q_{t_e})}\non\\
 & \le & \bar C\|\varphi\|_{L^4(Q_{t_e})},
\end{eqnarray}
where we applied the H\"older and Young inequalities. Thus we arrive
at $\|S_t\|_{L^\frac43(Q_{t_e})}\le\bar C$,   thus prove
(\ref{es2b}).
  And the proof of this lemma  is complete.

\medskip
For  the solution to the elliptic part of the system, i.e.
(\ref{eq1'}) -- (\ref{eq1a'}), we have
\begin{Lemma} \label{L2.3}
There hold for almost every  $t\in [0,t_e]$ that
\begin{eqnarray}
 & &\|u(t)\|_{W^{1,p}(\Omega)} + \|T(t)\|_{L^{ p}(\Omega)}\le \bar
 C,\label{est1}\\[0.2cm]
 & & \|u(t)\|_{H^{2}(\Omega)} + \|T(t)\|_{H^{1}(\Omega)}\le \bar C .
 \label{est1v}
\end{eqnarray}

\end{Lemma}
\noindent{\it Proof.} Using the estimate (\ref{est2}), we get $S(t)\in L^p(\Omega)$ for almost
 every $t\in [0,t_e]$ since the domain $\Omega$ is bounded. Recalling estimate (\ref{es2a}),
we obtain easily (\ref{est1}) -- (\ref{est1v}), by the regularity theory of
elliptic systems (or just using the formula \eq{Tformula} since our problem is one dimensional).
This completes the proof of the lemma.

\medskip
Now we differentiate (\ref{eq1appro}) once formally with respect to
$t$ and  recall the assumption
 on $b_t$, then use again the theory of the elliptic system to get
\begin{Lemma} \label{L2.5}
There hold for almost every  $t\in [0,t_e]$ that
\begin{eqnarray}
  \|T_t\|_{L^p(\Omega)}  \le  \bar  C\left(1+\|S_t\|_{L^p(\Omega)}\right),\
  \|T_t\|_{ L^\frac43(Q_{t_e}) } &\le& \bar  C,
 \label{est1t} \\
 T\in C([0,t_e];C^\alpha(\bar\Omega)),\quad {\rm and} \quad
 \|T\|_{C(\bar Q_{t_e})}&\le& \bar C.
 \label{est1a}
\end{eqnarray}

\end{Lemma}

To prove  the above lemma, we shall make use of the following lemma
which is of Aubin-Lions type, see, for instance, Lions~\cite{Lions},
and for the case $r=1$, see Simon~\cite{Simon87},
Roub\'icek~\cite{Roubicek}.
\begin{Lemma} \label{L2.6a}
Let $B_0,\ B,\ B_1$ be Banach spaces which satisfy
that $B_0,\ B_1$ are reflexive and that
$$
B_0\subset\subset B\subset B_1.
$$
Here, by $\subset\subset$ we denote the compact imbedding.
Define
$$
 W=\left\{f \mid f\in L^\infty(0,t_e;B_0), \quad
 \frac{df}{dt}\in L^r(0,t_e;B_1) \right\}
$$
with $t_e$ being a given positive number and $1<r<\infty$.

Then the embedding of $W$ in $C([0,t_e]; B)$ is compact.

\end{Lemma}

\noindent {\it Proof of Lemma 2.5.} We need only to prove (\ref{est1a}). From
(\ref{est1t}) and Lemma~\ref{L2.3}, we assert that
\begin{eqnarray}
T\in L^\infty(0,t_e;H^{1 }(\Omega)),\quad {\rm and} \quad T_t\in
L^\frac43(0,t_e;L^{\frac43}(\Omega)).
\end{eqnarray}
Thus we can choose
$$
 B_0=H^{1}(\Omega),\quad B=C^{\alpha}(\bar\Omega), \quad
 B_1=L^\frac43(\Omega),\quad r=\frac43,
$$
which meet the requirements of Lemma~\ref{L2.6a}, and $\alpha\in
(0,\frac12]$. Whence (\ref{est1a}) holds.  And  the proof of
Lemma~\ref{L2.5} is complete.

\medskip
Furthermore, for any fixed $\kappa$, we have
\begin{Lemma} \label{L2.7}
There hold  for any $t\in [0,t_e]$ that
\begin{eqnarray}
 \|S_t(t)\|^2 + \int_0^t\int_\Omega\left(|S_x|_\kappa+1\right)|S_{xt}|^2 dxd\tau  &\le& C,
 \label{approest3}  \\
 \|S_{xx}(t)\|&\le& C.
 \label{approest3a}
\end{eqnarray}

\end{Lemma}

\noindent {\it Proof.}  We prove firstly (\ref{approest3}). To this end, we
differentiate equation (\ref{eq2appro}) formally  with respect to $t$, then multiply
the resulting equation by $S_t$ and integrate it with respect to $x$ to get
\begin{eqnarray}
 \frac12\frac{d}{dt}\|S_t\|^2-c\,\nu\int_\Omega\left( |S_x|_\kappa S_{xx}\right)_tS_t
 \, dx  - c\int_\Omega \left( (T\cdot\bar\varepsilon-\hat\psi^\prime(S) )
 (|S_x|_\kappa - \kappa)\right)_tS_t\, dx = 0.
 \label{approest1b}
\end{eqnarray}
It is easy to see that
\begin{eqnarray}
 \left( |S_x|_\kappa S_{xx}\right)_t &=& \left(\int^{S_x} |\xi|_\kappa  d\xi \right)_{xt}, \\
 \int_\Omega\left(\int^{S_x} |\xi|_\kappa  d\xi\right)_{t}S_{xt}\, dx &=&
 \int_\Omega |S_x|_\kappa |S_{xt}|^2dx.
 \label{approest1a}
\end{eqnarray}
Thus \eq{approest1b} becomes
\begin{eqnarray}
 \frac12\frac{d}{dt}\|S_t\|^2
  + c\, \nu\int_\Omega  |S_x|_\kappa |S_{xt}|^2  dx - c\int_\Omega \left( (T\cdot\bar\varepsilon-\hat\psi^\prime(S))
 (|S_x|_\kappa-\kappa) \right)_tS_t\, dx   = 0.
 \label{approest1}
\end{eqnarray}
We now handle the last term on the left-hand side of (\ref{approest1})  as
\begin{eqnarray}
 && \left|c\int_\Omega \left( \left(T\cdot\bar\varepsilon-\hat\psi^\prime(S)\right)
 (|S_x|_\kappa - \kappa)\right)_tS_t\, dx\right|\non\\
 &\le & C\int_\Omega \left|T_t\cdot\bar\varepsilon-\hat\psi''(S)S_t\right|  (|S_x|_\kappa + 1)|S_t|dx
 +  C\int_\Omega \left|T\cdot\bar\varepsilon-\hat\psi'(S)\right|
 \, \left|(|S_x|_\kappa)'S_{xt}S_t \right|dx\non\\
 &\le & C\left(\|T_t\|+ \|S_t\|\right) (\|\,|S_x|_\kappa\|_{L^\infty(\Omega)} + 1) \|S_t\|
 +C\|S_t\|\,\|S_{xt}\|\non\\
 &\le & C\Big(\left(1+\|\,|S_x|_\kappa\|_{L^\infty(\Omega)}\right) \|S_t\|^2
 + (\|\,|S_x|_\kappa\|_{L^\infty(\Omega)}+ 1) \|T_t\| \, \|S_t\|  \Big)
 +\frac{\kappa}2 \|S_{xt}\|^2.
 \label{approest2}
\end{eqnarray}
Here we used the estimates $\|(T,\, S )\|_{L^\infty(Q_{t_e})}\le C,\ \|T_t\| \le C(1 + \|S_t\|),
\ \left|( |y|_\kappa)'\right|\le C$, and the H\"older, Young inequalities.
By the Sobolev imbedding theorem, we have
\begin{eqnarray}
 \|\,|S_x(t)|_\kappa\|_{L^\infty(\Omega)}\le C(\|S_x(t)\|_{L^\infty(\Omega)}+ 1) \le
 C(\|S_x(t)\|_{H^1(\Omega)} + 1).
 \label{Remark2.3}
\end{eqnarray}
Applying the estimate (\ref{est1t}), which is  valid for  $p=2$,
 combining (\ref{approest2}) and (\ref{approest1}), we arrive at
\begin{eqnarray}
 \frac{d}{dt}\|S_t(t)\|^2 \le C\left(1+\|S_x(t)\|_{H^1(\Omega)}\right) \|S_t(t)\|^2
 + C\left(\|S_x(t)\|_{H^1(\Omega)}^2 + 1\right).
 \label{approest2a}
\end{eqnarray}
We shall make use of the Gronwall inequality in the following form
\begin{eqnarray}
 y'(t)\le A(t)y(t)+B(t)\quad {\rm implies} \quad y(t)\le y(0) {\rm
 e}^{\int_0^t A(\tau)d\tau} + \int_0^t B(s){\rm e}^{\int_s^t
 A(\tau)d\tau}ds,
 \label{gronwall}
\end{eqnarray}
where $y,\, A,\, B$ are functions satisfying that
$ y(t)\ge 0$, $ A(t),\, B(t)$ are integrable over $[0,t_e]$. Defining
$$
 y(t)=\|S_t(t)\|^2, \quad A(t)= C\left(1+\|S_x\|_{H^1(\Omega)}\right),\quad
 B(t)=C\left(1+\|S_x\|_{H^1(\Omega)}^2\right) ,
$$
from the estimate $\|S_{x}\|_{H^1(Q_{t_e})}\le C$ which is a consequence of (\ref{es2a}) and the fact $|p|_\kappa\ge \kappa$,
it follows that the above defined $A(t),\, B(t)$ are integrable over $[0,t_e]$.
Thus we  can apply (\ref{gronwall}) to (\ref{approest2a}) and  obtain
$$
 \|S_t(t)\|^2\le C,
$$
whence
\begin{eqnarray}
 \|S_t(t)\|^2 + \int_0^t\int_\Omega |S_x|_\kappa |S_{xt}|^2 dx d\tau \le C.
\end{eqnarray}
From which we obtain easily (\ref{approest3}) since $\kappa$ at this moment
is a given number. Therefore we can use the equation to get easily (\ref{approest3a}).
Thus the proof of this lemma is complete.

\medskip
\noindent{\bf Remark 2.3.} {\it  Since we use the imbedding \eq{Remark2.3} which is valid only in one
dimensional case, thus this lemma is only true for this one dimensional problem.
}

\medskip
\noindent{\bf Remark 2.4.} {\it To derive \eq{approest3} rigorously, we employ the technique of
finite difference as, e.g., in \cite{Dafermos}. We assume that there exists a unique classical solution
$(u,T,S)$ to problem \eq{eq1appro} -- \eq{IDappro} such that
$$
 (u,T,S)\in C^{2,1}(\bar{Q}_{t_e}) \times C^{1,1}(\bar{Q}_{t_e})
 \times C^{2+\alpha,1+\alpha/2}(\bar{Q}_{t_e}), \quad S_{xt} \in L^2(Q_{t_e}).
$$
Define $S_h (t,x) = (S(t+h,x) - S(t,x))/h$ for any $h>0$. Then from \eq{eq2appro} we obtain
\begin{eqnarray}
 S_{ht} &=& \frac{c\, \nu}h\left(\int_{S_x(t,x)}^{S_x(t+h,x)} |\xi|_\kappa d\xi\right)_{x}
 + \left.\frac{c} h (T\cdot \bar\varepsilon-\hat \psi^\prime(S) )(|S_x|_\kappa - \kappa)\right|_{(t,x)}^{(t+h,x)},
 \label{justify1}
\end{eqnarray}
for any $[t\in [0,t_e-\delta]$, where $\delta$ is a fixed number such that $\delta\ge h$.
Here and hereafter, we use the notations $\big.f\big|_{(t,x)}^{(t+h,x)} = f (t+h,x) - f (t,x) $ and $f|_{(t,x)} = f (t,x) $.
Multiplying \eq{justify1} by $S_h$ and integrating the resulting equation with respect to $t,x$ over $Q_{t}$
yield
\begin{eqnarray}
 && \|S_{h }(t)\|^2  + c\, \nu\int_0^{t_e-\delta}\int_\Omega\frac1h \int_{S_x(t,x)}^{S_x(t+h,x)} |\xi|_\kappa d\xi S_{hx} dx dt \non\\
 &=& \|S_{h }(0)\|^2 + c\int_0^{t_e-\delta}\int_\Omega   (T_h\cdot \bar\varepsilon-\hat \psi^{\prime\prime}(S^*)S_h )
  (|S_x|_\kappa - \kappa)|_{(t+h,x)} S_h  dx dt \non\\
 && +\, c\int_0^{t_e-\delta}\int_\Omega (T\cdot \bar\varepsilon-\hat \psi^\prime(S) )|_{(t,x)}\,
  \frac{S_x(t+h,x) + S_x(t,x)}{|S_x(t+h,x)|_\kappa + |S_x(t,x)|_\kappa}   S_{xh} S_{h} dx dt.
 \label{justify2}
\end{eqnarray}
Here $S^*$ is a number between $S(t+h,x)$ and $S(t,x)$. Note that the second term on the left hand side of  \eq{justify2} is equal to
$$
 c\, \nu\int_0^{t_e-\delta}\int_\Omega |S_x(\eta,x) |_\kappa |S_{hx}|^2 dxdt \ge C\|S_{hx}\|^2_{L^2(Q_{t_e -\delta})},
$$
where $\eta\in [t,t+h]$ and we used $| p |_\kappa\ge \kappa$. By definition, we have
$$
 \frac{\left|S_x(t+h,x) + S_x(t,x) \right| }{|S_x(t+h,x)|_\kappa + |S_x(t,x)|_\kappa} \le 1.
$$
So the second  and third terms
 on the right hand side of \eq{justify2} are of lower orders
and can be estimated in a similar way  to \eq{approest2}. We thus arrive at
$$
 \|S_{hx}\|^2_{L^2(Q_{t_e-\delta} )}\le C.
$$
Further, we write
\begin{eqnarray}
 && \int_0^{t_e-\delta}\int_\Omega |S_x(\eta,x) |_\kappa |S_{hx}|^2 dx dt \non\\
 &=& \int_0^{t_e-\delta}\int_\Omega \Big( (|S_x(\eta,x) |_\kappa - |S_x(t,x) |_\kappa) + |S_x(t,x) |_\kappa \Big) |S_{hx}|^2 dx dt  .
 \label{justify3}
\end{eqnarray}
Invoking the H\"older continuity of $S_x$, applying the Fatou lemma for any fixed $\delta$, we take the limit as $ h\to 0$. Then letting $\delta \to 0$,
%to the non-negative sequence $\{(|S_x(t,x) |^\frac12_\kappa S_{hx})^2\}_h$, we obtain
%$$
% \int_0^{t_e-\delta}\int_\Omega  |S_x(t,x) |_\kappa  |S_{tx}|^2 dx dt\le \lim\inf_{h\to 0}I_2.
%$$
 we justify \eq{approest3} and omit  details.
}

\bigskip
We now turn back to prove  Theorem~\ref{T2.1}.

\noindent{\it Proof of Theorem~\ref{T2.1}.}  To complete the proof of the global existence of classical solution,
we need to prove that $\|S_x\|_{C^{\alpha/2,\alpha(\bar Q_{t_e}) }}\le C$.
To this end we make use of the estimates listed in  Lemmas~\ref{L2.2} -- ~\ref{L2.5} and Lemma~2.7.
To prove this, we invoke the following lemma see, e.g. \cite{Ladyzenskaya}
%
%%%%%%%%%%%%%%%%%%%%%%%%%%%%%%%%%%%%%%%%%%%%%%%%%%%%%%%%%%%%%%%%%%%%%%%%%%%%
%%%%%%%%%%%%%%%%%%%%%%%%%%%%%%%%%%%%%%%%%%%%%%%%%%%%%%%%%%%%%%%%%%%%%%%%%%%%
%
\begin{Lemma} \label{L2.8}
Let $f(t,x)$ be a function on $Q_{t_e}$ such that

i) $f$ is uniformly (with respect to $x$) H\"older continuous in $t$,
with exponent

$0<\alpha\le 1$, that is
$|f(t,x) - f(s,x)|\le C|t-s|^\alpha $,

and

ii) $f_x$ is uniformly (with respect to $t$) H\"older continuous in $x$,
with exponent

$0<\beta \le 1$, that is
$|f_x(t,x) - f_x(t,y)|\le C^\prime |y-x|^\beta $.

Then $f_x$ is uniformly H\"older continuous in $t$ with exponent $\gamma=\alpha\beta/(1+\beta)$,
such that
$$
 |f_x(t,x) - f_x(s,x)|\le C^{\prime\prime} |t-s|^\gamma ,\
 \forall x\in \bar\Omega, 0\le s\le t \le {t_e}.
$$
where $C^{\prime\prime}$ is a constant which may depend on $C, C^{\prime}$ and $\alpha, \beta$.

\end{Lemma}

By applying this lemma we assert that there exists a constant  $0<\alpha<1$ such that
$\|S_x\|_{C^{\alpha/2,\alpha }}\le C$. By the {\it a priori} estimate of the Schauder
type for parabolic equations,  we thus obtain that
$$
 \|S\|_{C^{1+\alpha/2,2+\alpha}}(\bar Q_{t_e})\le C,
$$
which ensures us to apply the
 Leray-Schauder fixed point theorem, and the proof of global existence
of classical solution is complete. Using
the technique of difference quotient with respect to $t$ to this classical solution see
e.g. \cite{Ladyzenskaya} we can prove (\ref{Sxt}). And the proof of Theorem~\ref{T2.1}
 is thus complete.

%%%%%%%%%%%%%%%%%%%%%%%%%%%%%%%%%%%%%%%%%%%%%%%%%%
%%%%%%%%%%%%%%% Section 3 %%%%%%%%%%%%%%%%%%%%%%%%
%
\section{Existence of weak solutions}

\subsection{Uniform a priori estimates}

This subsection is devoted to derivation of some uniform  {\it a priori }
estimates, which are independent of $\kappa\in (0,1]$,  for the approximate solution to
(\ref{eq1appro}) -- (\ref{IDappro}).  However these
estimates may depend on $\nu$, this thus makes it difficult to discuss the sharp
interface limit $\nu\to 0$. To investigate such a sharp interface limit,
we need new techniques.

We now denote the approximate solution by $(u^\kappa,T^\kappa,S^\kappa)$. Therefore we collect a priori estimates,
which have been established in Section~2 and are independent of $\kappa$.
\begin{Lemma} \label{L3.1}
There hold  for any $t\in [0,t_e]$ that
\begin{eqnarray}
 \|S_x^\kappa(t)\|^2 + \int_0^t\int_\Omega | S_x^\kappa|_\kappa |S_{xx}^\kappa|^2 dx d \tau &\le& \bar C,
 \label{est2a}\\
 \int_0^t\int_\Omega \left(\left(| S_x^\kappa|_\kappa |S_{xx}^\kappa|\right)^\frac43 +  |S^\kappa_t |^\frac43 \right)dx d \tau  &\le& \bar C
 \label{est2b},\\
 \int_0^t\|S^\kappa\|^2_{H^1(\Omega)}d\tau &\le& \bar C.
 \label{est2c}
\end{eqnarray}

\end{Lemma}

\medskip
\noindent{\bf Remark 3.1.} {\it From \eq{Remark31} we see that the constant $\bar C$ depends on $\nu$.
}
%

%%%%%%%%%%%%%%%%%%%%%%%%%%%%%%%%%%%%%%%%%%%%%%%%%%
%%%%%%%%%%%%%%% Section 4 %%%%%%%%%%%%%%%%%%%%%%%%
%
\subsection{Limits}

With the help of Lemma~2.6, applying the uniform a priori estimates established in Subsection~3.1, we shall investigate in this section
the limits, as $\kappa\to 0$, of the approximate solutions and complete the proof of
Theorem~\ref{T1.2}.
%We are going to make use of the {\it a priori }
%estimates in
%

\vskip0.2cm
\noindent{\it Proof of Theorem~\ref{T1.2}.} Firstly we
apply again Lemma~\ref{L2.6a} to show that the sequence of the
approximate solution $S^\kappa$ has a subsequence which converges
strongly. To this end, we choose
$$
 B_0=H^1(\Omega),\quad  B=C^\alpha(\bar\Omega),\quad  B_1=L^\frac43(\Omega),
$$
and
$$
 0<\alpha<\frac12,\quad p_1=\frac43.
$$
It is easy to see that such defined $B_0,\ B_1$ are reflexive.
Therefore we  apply Lemma~\ref{L2.6a} and conclude that the
sequence $\{S^\kappa\}_\kappa$ is a compact in
$C([0,t_e];C^\alpha(\bar\Omega))$. Thus we can
 select a subsequence of it,  and denote it by $\{S^{\kappa_n}\}_n$,
such that, as $n\to\infty$,
$$
 {\kappa_n}\to 0,
$$
and
\begin{eqnarray}
 S^{\kappa_n}\to S,\ {\rm in}\ C([0,t_e];C^\alpha(\bar\Omega)),
 \label{aecong1}
\end{eqnarray}
from which we obtain that
\begin{eqnarray}
 \|S^{\kappa_n} - S\|_{C([0,t_e]\times \bar\Omega)}\to 0.
 \label{aecong}
\end{eqnarray}
 On the other hand, by Lemma~2.5, we assert that
 there exists a subsequence of $T^\kappa$ such that
\begin{equation}
 T^{\kappa_n} \to T\ {\rm in}\ C^\alpha(\bar Q_{t_e}),
 \label{convergence0}
\end{equation}
from this, (\ref{aecong}) and \eq{uformula},  we obtain consequently that
\begin{equation}
 (u^{\kappa_n} ,T^{\kappa_n})\to (u,T), \ {\rm uniformly\ in}\ C^\alpha(\bar Q_{t_e}),
 \label{convergence0a}
\end{equation}
as $n\to \infty$.

\bigskip
In what follows, we are going to prove that the limit function
 $(u,T,S)$ is just a weak solution to problem (\ref{eq1a}) -- (\ref{eq4a})
 in the sense of Definition~\ref{D1.1}.
 It is not difficult to show that (\ref{eq1a}) and (\ref{eq1aa})
are satisfied by the linearity of those two equations and
 the uniform convergence of $u^\kappa,\, T^\kappa$. The remaining part of this section is
 to prove that (\ref{eq2a}) is satisfied.

We shall make use of the theorem on the stability of viscosity solutions, see e.g. \cite{Giga}.
%%%%%%%%%%%%%%%%%%%%%%%
%%%%%%%%
\begin{theorem}[Stability of viscosity solutions] Assume that $F_n$ converges to $F$
locally uniformly (as $n\to\infty$) in the domain of definition of $F$. Assume
that $v_n$ is a viscosity solution to
$$
 (v_n)_t + F_n(t,x, v_n,\nabla v_n,\nabla^2 v_n) \le 0\ (resp. \ge 0)\ in\ Q_{t_e},
$$
and that $v_n$ converges to $v$ locally uniformly in $Q_{t_e}$ as $n\to\infty$.
Then $v$ is a viscosity solution to
$$
 v_t + F (t,x, v ,\nabla v,\nabla^2 v) \le 0\ (resp. \ge 0)\ in\ Q_{t_e}.
$$

\end{theorem}

\medskip
To apply this theorem, we define $v_n = S^{\kappa_n}$ and
$$
 F_n (t,x,p,q,r) = H_{T^{\kappa_n}}(t,x,p,q,r) = c\,\left(T^{\kappa_n}(t,x)\cdot\bar\varepsilon -\hat\psi'(p)+\nu r\right)(|q|_{\kappa_n} - {\kappa_n}).
$$

Invoking \eq{convergence0a} and \eq{aecong} we conclude that
\\[0.12cm]
i) $S^{\kappa_n}$ converges to $S$ locally uniformly in any compact subset in $ Q_{t_e} $.
\\[0.12cm]
ii) $H_{T^{\kappa_n}} (t,x,p,q,r)$ converges to $H_T$ locally uniformly in any compact subset in $(0,t_e)\times\Omega\times\cR\times \cR\times\cR $.
 \\[0.12cm]
iii) Since $S^{\kappa_n}$ is a classical solution to equation (\ref{eq2appro}) when  $T^{\kappa_n}$ is regarded temporarily fixed, $S^{\kappa_n}$ is
also a viscosity solution to (\ref{eq2appro}).

Therefore, we can apply  Theorem~4.1 and conclude that the limit $S$ is a viscosity solution to
$S_t = H_T(t,x, S_x, S_{xx})$. Hence, recalling the properties $S^*(t,x)=S_*(t,x)=S(t,x)$ and $(H_T)^*( t,x,p,q,r)=(H_T)_*( t,x,p,q,r)=  H_T ( t,x,p,q,r)$, we assert that  $(u,T,S)$ is a weak solution to problem (\ref{eq1a}) -- (\ref{eq4a})
 in the sense of Definition~\ref{D1.1}, and the  proof of Theorem~\ref{T1.2} is thus complete.

\bigskip
\section{Appendix}

Since our model is quite new, we briefly sketch, for the sake of readers' convenience,
 the physical background and the derivation of the diffusive interface model \eq{eq1} -- \eq{eq3}
 from a sharp interface model. We also refer the reader to \cite{Alber06,Alber04,Alber07}.
 Our model differs from the Allen-Cahn  model by a gradient term. The main reason is: In the Allen-Cahn model, the driving force for
the motion of interface is  the mean curvature, while the motion of interface considered in this
article is driven by configurational forces, see e.g. \cite{Gurtin,Maugin}.

 Material phases are characterized by the structure of the crystal lattice, in which
the atoms are arranged. An interface between different material phases moves if
the crystal lattice in front of the interface is transformed from one structure to the
other. Often phase transformations are triggered by diffusion processes. A well-known
model for diffusion dominated transformations is the Allen-Cahn equation when the order parameter
is not conserved (or the Cahn-Hilliard equation if the order parameter is conserved).
We derive our model \eq{eq1} -- \eq{eq3} from a sharp interface model for diffusionless transformations, also called
martensitic transformations, see e.g. \cite[p. 162]{Horn}. This sharp interface model is an
initial-boundary value problem for the unknown functions $u,\ T$
and for the unknown interface $\Gamma(t)\subset \Omega$ between two material phases, which is a free boundary.
It consists of \eq{eq1} -- \eq{eq02} and the interface conditions
\begin{eqnarray}
 V(t, x)[S](t, x) &=& c\left( -\langle T\rangle(t,x)\cdot\bar\varepsilon [S] (t, x)
 + [\hat\psi(S)](t, x)\right),
 \label{A1}\\
 \ [u](t, x) &=& 0,\ [T](t, x)n(t, x) = 0,
 \label{A2}
\end{eqnarray}
which must hold for $x\in\Gamma (t)$, and of a Dirichlet boundary condition for $u$  the
initial condition  \eq{eq3}. We use the notation $[f] = f_+ - f_-$ and $\langle f \rangle
 = \frac12 (f_+ + f_-)$, where $f_+,\ f_-$ are the limit values of the function $f$ on both sides of $\Gamma(t)$. Moreover,
$V(t, x)\in\cR^3$ denotes the normal speed of the interface $\Gamma(t)$, which is measured as
positive in the direction for which $[S](t, x)$ is positive. Here $c$ is a positive constant.
Equation \eq{A1}, a constitutive equation, determines the normal speed $V$ of the
phase interface as a function of the term $-\langle T\rangle \cdot\bar\varepsilon [S] + [\hat\psi(S)]$. Some computations
show that this term is equal to the expression $n \cdot [E]n$ with the Eshelby tensor $E$ (an energy-momentum
tensor, see \cite[p753-p767]{Eshelby}) and the normal vector $n$ to  $\Gamma(t)$ (cf. \cite{Alber04}) and thus is a configurational force. We
assume that $V$ depends linearly on the configurational force, which is the most simple
constitutive assumption. Thus, in this model the evolution of the phase interface is
driven by the configurational force along the interface, an assumption appropriate for
martensitic transformations.

Though configurational forces were introduced in the first half of the last century,
it was clearly stated for the first time in \cite{AK90} that \eq{eq1}, \eq{eq02}, \eq{A1}, \eq{A2} form a closed
initial-boundary value problem. Applications of this model can be found, for example,
in \cite{Buratti,Mueller,Socrate}, where equilibrium configurations for materials with phase transitions
are determined, and in \cite{James}, where the evolution of phase interfaces in ferroelectric
materials is modeled. In a sense, this free initial-boundary value problem from solid
mechanics is comparable to the Stefan problem in fluid mechanics.

The initial-boundary value problem \eq{eq1} -- \eq{eq3} can be considered to be a regularization
of this sharp interface model, which could be used to prove existence of
solutions of the sharp interface model, and it can also be considered to be a diffusive
interface model for martensitic phase transitions, which is useful by itself and avoids
some disadvantages of the model with sharp interfaces. We are interested in both
aspects.

The derivation of \eq{eq1} -- \eq{eq3} given in \cite{Alber00,Alber04}  uses a rigorous method. To make
the model plausible, we derive the model here in a different, short, but formal way.
To this end we replace the phase interface $\Gamma(t)$, across which the order parameter
jumps from $0$ to $1$, by finitely many interfaces parallel to the original interface, and
consider a new order parameter, again denoted by $S$, with small jumps across these
interfaces, such that the sum of the jumps is equal to $1$. We assume that the new order
parameter satisfies \eq{A1} and \eq{A2} along all interfaces. If we increase the number of
interfaces and decrease the jump height, the new order parameter will converge to a
continuous or even differentiable order parameter, for which the normal speed of the
level manifolds is equal to the limit of the normal speed of the interfaces. For this
limit speed we obtain from \eq{A1}
\begin{eqnarray}
 V(t, x) &=& c \lim_{[S]\to 0}(\langle T\rangle \cdot\bar\varepsilon + \hat\psi'(S^*))
 = c\, (-T \cdot\bar\varepsilon   +  \hat\psi' (S)) = c\, \psi_S(\varepsilon(\nabla_xu), S).
\end{eqnarray}
The limit order parameter thus satisfies the Hamilton-Jacobi transport equation
\begin{eqnarray}
 S_t = - c\, \psi_S(\varepsilon(\nabla_xu), S) |\nabla_x S|,
 \label{A3}
\end{eqnarray}
since the level manifolds of solutions of  equation \eq{A3} have this normal speed.
The idea suggests itself to approximate the solution of the sharp interface model
by smooth solutions $(u, T, S)$ of the system \eq{eq1}, \eq{eq2}, \eq{A3}. Yet, examples in one
space dimension show that in general the function $S$ in such a smooth solution develops
a jump after finite time. The reason for this is that the function $\hat\psi'$
appearing in $ \psi_S$
is not monotone, since $\hat\psi$ is a double well potential. After $S$ has developed a jump,
\eq{A3} can no longer be used to govern the evolution of $S$. To avoid this problem and
to force solutions to stay smooth, \eq{A3} has been replaced by \eq{eq2}, which contains
the regularizing term $\nu |\nabla_x S|\Delta_x S$   with the small positive parameter $\nu$. This yields
the model \eq{eq1} -- \eq{eq3}.

The choice of this special regularizing term follows from the second law of thermodynamics,
which every model must satisfy. This law requires that there exist a
free energy $\psi$ and a flux $q$ such that $
\frac{\partial}{\partial t} \psi + {\rm div}_x\, q \le  b\cdot u_t$ holds; cf. \cite{Alt}. If we choose
 a free-energy and a flux as \eq{freeenergy} and \eq{flux},
it follows by a short computation for solutions $(u, T, S)$ of \eq{eq1}, \eq{eq2} that
$$
\frac{\partial}{\partial t} \psi - {\rm div}_x \left( Tu_t + \nu S_t \nabla_x S \right) - b\cdot u_t
= (\psi_S(\varepsilon, S) - \nu\Delta_x S )S_t.
$$
Inserting \eq{eq3} into this equation shows that the right-hand side is non-positive, whence
the second law is fulfilled. However this would {\it not} be true by using, as in the theory of
conservation laws, the standard regularization  (i.e. adding an artificial viscosity term)
$S_t = -c\, \psi_S (\varepsilon(\nabla_xu), S) |\nabla_x S| +\nu\Delta_x S$ of \eq{A3}.

\bigskip
\bigskip
\noindent{\bf Acknowledgement.}
\vskip0.12cm The author would like
to express his sincere thanks to Prof. H.-D. Alber for valuable
discussions.
%, and to the referee(s) for valuable opinions.
This work has been partly supported by Grant
MTM2008-03541 of the Ministerio de Educac\'ion y Ciencia of Spain, and by
Project PI2010-04 of the Basque Government.

%
%%%%%%%%%%%%%%%%%%%%%%%%%%%%%%%%%%%%%%%%%%%%%%%%%%%%%%%%%%%%%%%

\bigskip

\bibliographystyle{plainnat}

\end{document}